\documentclass[12pt]{article}
\usepackage{amsfonts}

\usepackage{graphicx}
\usepackage{amsmath}
\usepackage[USenglish]{babel}
\usepackage{a4}

\newtheorem{Prop}{Proposition}[section]
\newtheorem{Theo}[Prop]{Theorem}
\newtheorem{Rem}[Prop]{Remark}

\newtheorem{Def}[Prop]{Definition}
\input tcilatex

\begin{document}

\title{A non-stationary model for catalytic converters with cylindrical geometry}
\author{Jean-David\ HOERNEL \\
Laboratoire de Math\'{e}matiques et Applications\\
Universit\'{e} de Haute-Alsace, 4 rue des Fr\`{e}res Lumi\`{e}re\\
F-68093\ MULHOUSE Cedex FRANCE\\
Email : j-d.hoernel@wanadoo.fr}
\date{}
\maketitle

\begin{abstract}
We prove some existence and uniqueness results and some qualitative
properties for the solution of a system modelling the catalytic conversion
in a cylinder. This model couples parabolic partial differential equations
posed in a cylindrical domain and on its boundary.
\end{abstract}

\section{Introduction}

A gas containing $N-1$ different chemical species is flowing through a
cylindrical passage with a parabolic speed profile. Chemical species are
diffusing in the cylinder and are reacting only on the boundary of the
cylinder.

\medskip We investigate the existence, uniqueness and qualitative properties
of a non-stationary system of partial differential equations describing the
evolution of the concentrations of the $N-1$ chemical species and of the
temperature, both with respect to the time variable and along the cylinder.
This model of catalytic converters starts with the contribution of Ryan,
Becke and Zygourakis \cite{Ryan}.\ However, we have added an axial diffusion
term on the boundary in the present model, see \cite{these} for a more
detailed description of this model.

\medskip Due to its internal symetry, the cylinder may be reduced to the
domain $\Omega =\left[ 0,1\right[ \times \left] 0,1\right[ $, the boundary
of which is $\Sigma =\left\{ 1\right\} \times \left] 0,1\right[ $. The
concentrations (resp. the temperature) inside the cylinder $\Omega $ are
named $C_{if}$, $i=1,...,N-1$ (resp.$\;C_{Nf}$). The concentrations (resp.
the temperature) on the boundary $\Sigma $ are named $C_{is}$ (resp. $C_{Ns}$%
). The problem is written in a normalized way as:
\begin{equation}
\left\{
\begin{array}{rcl}
\left( 1-r^{2}\right) \dfrac{\partial C_{if}}{\partial z}\left( r,z,t\right)
& = & \dfrac{\beta _{if}}{r}\dfrac{\partial }{\partial r}\left( r\dfrac{%
\partial C_{if}}{\partial r}\right) \left( r,z,t\right) , \\
\dfrac{\partial C_{is}}{\partial t}\left( z,t\right) & = & -\gamma _{is}%
\dfrac{\partial C_{if}}{\partial r}(1,z,t)+\delta _{i}\mathbf{r}_{i}\left(
C_{1s}^{+},\ldots ,C_{Ns}^{+}\right) \left( z,t\right) \\
&  & \qquad +\theta _{is}\dfrac{\partial ^{2}C_{is}}{\partial z^{2}}\left(
z,t\right) ,
\end{array}
\right.  \label{sys1}
\end{equation}
for $i\in \left\{ 1,\ldots ,N\right\} $ and with $C_{is}^{+}\left(
z,t\right) =\sup \left( C_{is},0\right) $.\ The initial and boundary
conditions are:
\begin{equation}
\left\{
\begin{array}{rclrrcl}
C_{if}\left( r,0,t\right) & = & C_{i0}\left( r\right) , &  & \dfrac{\partial
C_{if}}{\partial r}\left( 0,z,t\right) & = & 0, \\
C_{if}\left( 1,z,t\right) & = & C_{is}\left( z,t\right) , &  & C_{is}\left(
z,0\right) & = & C_{is0}\left( z\right) , \\
\theta _{is}\dfrac{\partial C_{is}}{\partial z}\left( 1,t\right) & = & 0, &
& \theta _{is}\dfrac{\partial C_{is}}{\partial z}\left( 0,t\right) & = & 0.
\end{array}
\right.  \label{cond1}
\end{equation}

The functions $\mathbf{r}_{i}$, $i\in \left\{ 1,\ldots ,N\right\} $, are
supposed to be Lipschitz continuous:
\begin{equation*}
\left| \mathbf{r}_{i}\left( C_{1s}^{1},\ldots ,C_{Ns}^{1}\right) -\mathbf{r}%
_{i}\left( C_{1s}^{2},\ldots ,C_{Ns}^{2}\right) \right| \leq
k_{i}\sum\limits_{h=1}^{N}\left| C_{hs}^{1}-C_{hs}^{2}\right|
\end{equation*}
and verify the following hypotheses:

\begin{itemize}
\item[(H1)]  $\forall x\in \left( \mathbb{R}^{+}\right) ^{N}:\mathbf{r}%
_{i}\left( x_{1},\ldots ,x_{N}\right) \geq 0.$

\item[(H2)]  If at least one of the $x_{i}$, $1\leq i\leq N$, is equal to $0$%
, then $\mathbf{r}_{i}\left( x_{1},\ldots ,0,\ldots ,x_{N}\right) =0$.

\item[(H3)]  For every $\left( x,y\right) \in \left( \mathbb{R}^{+}\right)
^{N}\times \left( \mathbb{R}^{+}\right) ^{N}$:
\begin{equation*}
-\sum\limits_{i=1}^{N}\delta _{i}\dfrac{\beta _{if}}{\gamma _{is}}\left(
\mathbf{r}_{i}\left( x_{1},\ldots ,x_{N}\right) -\mathbf{r}_{i}\left(
y_{1},\ldots ,y_{N}\right) \right) \left( x_{i}-y_{i}\right) \geq 0.
\end{equation*}
\end{itemize}

\begin{Rem}
We observe that:

\begin{enumerate}
\item  Hypotheses $\left( H2\right) $ and $\left( H3\right) $ imply:
\begin{equation*}
\forall x\in \left( \mathbb{R}^{+}\right) ^{N}:-\sum\limits_{i=1}^{N}\delta
_{i}\dfrac{\beta _{if}}{\gamma _{is}}\mathbf{r}_{i}\left( x_{1},\ldots
,x_{N}\right) x_{i}\geq 0.
\end{equation*}

\item  We have:
\begin{equation*}
C_{if}\left( 1,z,t\right) =C_{is}\left( z,t\right) \text{ ; }C_{if}\left(
r,0,t\right) =C_{i0}\left( r\right) \text{ ; }C_{is}\left( z,0\right)
=C_{is0}\left( z\right) .
\end{equation*}

In order to ensure the continuity of the concentrations in $\Omega $ and on
the boundary $\Sigma $ at $z=0$, we must have at $t=0$, $C_{if}\left(
1,0,0\right) =C_{is}\left( 0,0\right) $, which implies : $C_{i0}\left(
1\right) =C_{is0}\left( 0\right) $ ; $C_{is}\left( 0,t\right) =C_{i0}\left(
1\right) $.
\end{enumerate}
\end{Rem}

\section{\textbf{Existence of the solution}}

We establish the existence of the solution using the diagram:
\begin{equation*}
\begin{array}{cccccc}
&  &  & \Gamma &  &  \\
& \tilde{C}_{if} &  & \rightarrow &  & C_{if} \\
\Phi & \downarrow &  & \nearrow &  &  \\
&  &  & \Psi &  &  \\
& C_{is} &  &  &  &
\end{array}
\end{equation*}

Indeed, the proof of the existence is decomposed in two steps:

\begin{enumerate}
\item  Existence of a solution in $\Omega $ (given the solution on the
boundary $\Sigma $) and existence of a solution on the boundary $\Sigma $
(given the solution in $\Omega $);

\item  $\Gamma =\Psi \circ \Phi $ is a contraction in some appropriate
functional space.
\end{enumerate}

\subsection{\textbf{Preliminary results}}

We have the following continuous embeddings (see \cite[p. 103]{Kuf}):
\begin{equation*}
L^{2}\left( \Omega \right) \underset{\rightarrow }{\subset }L_{r}^{2}\left(
\Omega \right) \underset{\rightarrow }{\subset }L_{r(1-r^{2})}^{2}\left(
\Omega \right) \text{ ; }W^{1,2}\left( \Omega \right) \underset{\rightarrow
}{\subset }W_{r}^{1,2}\left( \Omega \right) ;W^{1,2}\left( \Omega \right)
\underset{\rightarrow }{\subset }W_{r(1-r^{2})}^{1,2}\left( \Omega \right) .
\end{equation*}

Because $\frac{\partial C_{if}}{\partial r}(1,z,t)=\frac{1}{\beta _{if}}%
\int\nolimits_{0}^{1}\frac{\partial C_{if}}{\partial z}r\left(
1-r^{2}\right) dr$, we can write the problem as:
\begin{equation}
\left\{
\begin{array}{rll}
\left( 1-r^{2}\right) \dfrac{\partial C_{if}}{\partial z} & = & \beta _{if}%
\dfrac{1}{r}\dfrac{\partial }{\partial r}\left( r\dfrac{\partial C_{if}}{%
\partial r}\right) , \\
\dfrac{\partial C_{is}}{\partial t}-\theta _{is}\dfrac{\partial ^{2}C_{is}}{%
\partial z^{2}} & = & -\dfrac{\gamma _{is}}{\beta _{if}}\int%
\nolimits_{0}^{1}\dfrac{\partial C_{if}}{\partial z}r\left( 1-r^{2}\right)
dr+\delta _{i}\mathbf{r}_{i}(C_{1s}^{+},\ldots ,C_{Ns}^{+}),
\end{array}
\right.  \label{SR'}
\end{equation}
with $i\in \{1,\ldots ,N\}$, and the initial or boundary conditions:
\begin{equation*}
\left\{
\begin{array}{rlllrll}
C_{if}\left( r,0,t\right) & = & C_{i0}\left( r\right) , &  & \dfrac{\partial
C_{if}}{\partial r}\left( 0,z,t\right) & = & 0, \\
C_{if}\left( 1,z,t\right) & = & C_{is}\left( z,t\right) , &  & C_{is}(z,0) &
= & C_{is0}\left( z\right) , \\
\theta _{is}\dfrac{\partial C_{is}}{\partial z}\left( 1,t\right) & = & 0, &
& \theta _{is}\dfrac{\partial C_{is}}{\partial z}\left( 0,t\right) & = & 0.
\end{array}
\right.
\end{equation*}

\subsection{\textbf{Existence in the cylinder}}

Assuming that the $C_{is}$ are known on the boundary and performing a change
of function in order to have homogeneous boundary conditions at $r=1$, we
obtain the following problem in which we omitted the time variable $t$:
\begin{equation}
\left\{
\begin{array}{rll}
\left( 1-r^{2}\right) \dfrac{\partial u_{f}}{\partial z}-\beta _{f}\dfrac{1}{%
r}\dfrac{\partial }{\partial r}\left( r\dfrac{\partial u_{f}}{\partial r}%
\right) & = & -\dfrac{\partial u_{s}}{\partial z}, \\
u_{f}\left( r,0\right) & = & u_{0}\left( r\right) , \\
u_{f}\left( 1,z\right) & = & 0, \\
\dfrac{\partial u_{f}}{\partial r}\left( 0,z\right) & = & 0,
\end{array}
\right.  \label{SH}
\end{equation}
with $i\in \left\{ 1,\ldots ,N\right\} $ and
\begin{equation*}
\begin{array}{rrllrll}
w_{f} & = & ^{t}\left( C_{1f},\ldots ,C_{Nf}\right) , &  & u_{s} & = &
^{t}\left( C_{1s},\ldots ,C_{Ns}\right) , \\
u_{f} & = & w_{f}-u_{s}, &  & \beta _{f} & = & diag\left( \beta _{1f},\cdots
,\beta _{Nf}\right) .
\end{array}
\end{equation*}

We set :
\begin{equation*}
\begin{array}{lll}
W & = & \left\{ u\in \left( L_{r}^{2}(0,1)\right) ^{N}\mid \dfrac{\partial u%
}{\partial r}\in \left( L_{r}^{2}(0,1)\right) ^{N}\right\} , \\
W_{0} & = & \left\{ u\in \left( L_{r}^{2}(0,1)\right) ^{N}\mid \dfrac{%
\partial u}{\partial r}\in \left( L_{r}^{2}(0,1)\right) ^{N},u(1)=0\right\}
\end{array}
\end{equation*}
and let $W_{0}^{\prime }$ be the dual space of $W_{0}.$

\begin{Def}
\label{def_fluide}Assume that $u_{0}\in \left( L_{r}^{2}\left( 0,1\right)
\right) ^{N}$ and $\frac{\partial u_{s}}{\partial z}\left( z\right) \in
\left( L^{2}\left( 0,1\right) \right) ^{N}$. A function $u_{f}$ is a \textbf{%
weak solution} of (\ref{SH}) if and only if $u_{f}\in L^{2}(0,1;W_{0})$, $%
\frac{\partial u_{f}}{\partial z}\in L^{2}(0,1;W_{0}^{\prime })$ and if for
every $\varphi \in L^{2}(0,1;W_{0})$, we have:
\begin{equation*}
\begin{array}{l}
\int\nolimits_{0}^{1}\int\nolimits_{0}^{1}\left\langle \dfrac{\partial
u_{f}}{\partial z},\varphi \right\rangle
r(1-r^{2})drdz+\int\nolimits_{0}^{1}\int\nolimits_{0}^{1}\left\langle
\beta _{f}r\dfrac{\partial u_{f}}{\partial r},\dfrac{\partial \varphi }{%
\partial r}\right\rangle drdz \\
\qquad =-\int\nolimits_{0}^{1}\int\nolimits_{0}^{1}\left\langle \dfrac{%
\partial u_{s}}{\partial z},\varphi \right\rangle r(1-r^{2})drdz.
\end{array}
\end{equation*}
\end{Def}

\begin{Prop}
Let $u_{0}$ and $\frac{\partial u_{s}}{\partial z}$ be as in the precedent
definition. Then there exists at least one \textbf{weak solution} of (\ref
{SH}). This solution verifies:
\begin{equation}
\begin{array}{rll}
\left\| u_{f}\right\| _{L^{\infty }\left( 0,1;\left(
L_{r(1-r^{2})}^{2}\left( 0,1\right) \right) ^{N}\right) } & \leq & C, \\
\int\nolimits_{0}^{1}\left\| \dfrac{\partial u_{f}}{\partial r}(z)\right\|
_{\left( L_{r}^{2}\left( 0,1\right) \right) ^{N}}^{2}dz & \leq & C, \\
\left\| \dfrac{\partial u_{f}}{\partial z}\right\| _{L^{2}\left(
0,1;W_{0}^{\prime }\right) } & \leq & C,
\end{array}
\label{est1}
\end{equation}
where $C$ is a positive constant which only depends on the data of the
problem.
\end{Prop}

\begin{proof}
We use a Galerkin approximation of $u_{f}$ for which we establish the three
above estimates, and pass to the limit in order to prove the Proposition,
see \cite{these} for the details.
\end{proof}

\begin{Prop}
Under the conditions given in the Definition \ref{def_fluide}, the solution
of (\ref{SH}) is such that:
\begin{equation*}
u_{f}\in L^{2}(0,1;W_{0})\cap C\left( 0,1;\left( L_{r}^{2}\left( 0,1\right)
\right) ^{N}\right) .
\end{equation*}
\end{Prop}

\begin{proof}
From the inclusion $W_{0}\subset \left( L_{r}^{2}\left( 0,1\right) \right)
^{N}\subset W_{0}^{\prime }$, and the fact that $u_{f}\in L^{2}(0,1;W_{0})$,
$\frac{\partial u_{f}}{\partial z}\in L^{2}(0,1;W_{0}^{\prime })$, we deduce
the result using Proposition 23.23 of \cite[p. 422]{Zeidler}.
\end{proof}

\begin{Rem}
In order to take into account the time variable $t$, all the above
expressions $v\left( .\right) \in H$ have to be understood as $v\left(
.,t\right) \in L^{2}\left( 0,T;H\right) $.
\end{Rem}

\subsection{\textbf{Existence on the boundary}}

Assuming that $\frac{\partial u_{f}}{\partial z}$ is known, we have:
\begin{equation}
\left\{
\begin{array}{rll}
\dfrac{\partial u_{s}}{\partial t}-\theta _{s}\dfrac{\partial ^{2}u_{s}}{%
\partial z^{2}} & = & \delta \mathbf{r}(u_{s}^{+})-\Gamma
_{f}\int\nolimits_{0}^{1}\dfrac{\partial u_{f}}{\partial z}r\left(
1-r^{2}\right) dr, \\
u_{s}\left( z,0\right)  & = & u_{s0}\left( z\right) , \\
\theta _{s}\dfrac{\partial u_{s}}{\partial z}\left( 0,t\right)  & = & 0, \\
\theta _{s}\dfrac{\partial u_{s}}{\partial z}\left( 1,t\right)  & = & 0,
\end{array}
\right.   \label{sysbord1}
\end{equation}
with:
\begin{equation*}
\begin{array}{rrllrll}
\mathbf{r} & = & ^{t}\left( \mathbf{r}_{1,}\ldots ,\mathbf{r}_{N}\right) , &
& \Gamma _{f} & = & diag\left( \dfrac{\gamma _{1s}}{\beta _{1f}},\cdots ,%
\dfrac{\gamma _{Ns}}{\beta _{Nf}}\right) , \\
\delta  & = & \left( \delta _{1},\cdots ,\delta _{N}\right) , &  & \theta
_{s} & = & \left( \theta _{1s},\cdots ,\theta _{Ns}\right) .
\end{array}
\end{equation*}

Let
\begin{equation*}
H^{1}(0,1)=\left\{ u\in \left( L^{2}\left( 0,1\right) \right) ^{N}\mid
\dfrac{\partial u}{\partial z}\in \left( L^{2}\left( 0,1\right) \right)
^{N}\right\}
\end{equation*}
and let $H^{\ast }$ be the dual space of $H^{1}$.

\begin{Def}
\label{def_paroi}Suppose that $\frac{\partial u_{f}}{\partial z}$ belongs to
$L^{2}\left( 0,T;L^{2}\left( 0,1;\left( L_{r}^{2}\left( 0,1\right) \right)
^{N}\right) \right) $ and that $u_{s0}$ belongs to $\left( L^{2}\left(
0,1\right) \right) ^{N}$. A function $u_{s}$ is a \textbf{weak solution} of (%
\ref{sysbord1}) if and only if $u_{s}\in L^{2}\left( 0,T;H^{1}\left(
0,1\right) \right) ,$ $\dfrac{\partial u_{s}}{\partial t}\in L^{2}\left(
0,T;H^{\ast }\left( 0,1\right) \right) $, satisfies $u_{s}\left( z,0\right)
=u_{s0}\left( z\right) $, and if for all $\varphi \in L^{2}\left(
0,T;H^{1}\left( 0,1\right) \right) $, we have:
\begin{equation*}
\begin{array}{l}
\int\nolimits_{0}^{T}\int\nolimits_{0}^{1}\left\langle \dfrac{\partial
u_{s}}{\partial t},\varphi \right\rangle
dzdt+\int\nolimits_{0}^{T}\int\nolimits_{0}^{1}\left\langle \theta _{s}%
\dfrac{\partial u_{s}}{\partial z},\dfrac{\partial \varphi }{\partial z}%
\right\rangle dzdt \\
\qquad =\int\nolimits_{0}^{T}\int\nolimits_{0}^{1}\left\langle \delta
\mathbf{r}(u_{s}^{+}),\varphi \right\rangle
dzdt-\int\nolimits_{0}^{T}\int\nolimits_{0}^{1}\int\nolimits_{0}^{1}\left%
\langle \Gamma _{f}\dfrac{\partial u_{f}}{\partial z},\varphi \right\rangle
r\left( 1-r^{2}\right) drdzdt.
\end{array}
\end{equation*}
\end{Def}

\begin{Prop}
Let $u_{f}$ and $u_{s0}$ as in the preceding definition. Then there exists
at least one \textbf{weak solution} $u_{s}$ of (\ref{sysbord1})$.$
\end{Prop}

\begin{proof}
We prove the existence of a solution for the linearized weak formulation of
the problem due to Theorem 2.2 of \cite[p. 286]{Mal} and then we use some
fixed point argument.
\end{proof}

\begin{Prop}
Under the conditions given in the Definition \ref{def_paroi}, the solution
of (\ref{sysbord1}) is such that:
\begin{equation*}
u_{s}\in L^{2}\left( 0,T;H^{1}\left( 0,1\right) \right) \cap C\left(
0,T;\left( L^{2}\left( 0,1\right) \right) ^{N}\right) .
\end{equation*}
\end{Prop}

\begin{proof}
From the embeddings $H^{1}\left( 0,1\right) \subset \left( L^{2}\left(
0,1\right) \right) ^{N}\subset H^{\ast }\left( 0,1\right) $ and the fact
that $u_{s}\in L^{2}\left( 0,T;H^{1}\left( 0,1\right) \right) $, $\frac{%
\partial u_{s}}{\partial t}\in L^{2}\left( 0,T;H^{\ast }\left( 0,1\right)
\right) $, we deduce the result using Proposition 23.23 of \cite[p. 422]
{Zeidler}.
\end{proof}

\subsubsection{$\Gamma =\Psi \circ \Phi $ is a c\textbf{ontraction}}

Consider the mapping $\Gamma :\tilde{U}_{f}\mapsto U_{s}\rightarrow U_{f},$\
and let $\mu =\sup_{i}\left( \gamma _{is}/\beta _{if}\right)
^{1/2}/\inf_{i}\theta _{is}$.

\begin{Prop}
\begin{enumerate}
\item  The mapping $\Phi :\tilde{U}_{f}\rightarrow U_{s}$ is such that if we
suppose $\mu \alpha ^{2}<4$, then:
\begin{equation}
\int\nolimits_{0}^{\tau }\int\nolimits_{0}^{1}\left\| \dfrac{\partial U_{s}%
}{\partial z}\right\| ^{2}dzdt\leq \dfrac{4\mu }{\alpha ^{2}\left( 4-\mu
\alpha ^{2}\right) }\underset{s\in \left[ 0,1\right] }{\sup }%
\int\nolimits_{0}^{\tau }\int\nolimits_{0}^{1}\left\| \tilde{U}%
_{f}\right\| ^{2}(r,s,t)r\left( 1-r^{2}\right) drdt.  \label{phi}
\end{equation}

\item  The mapping $\Psi :U_{s}\rightarrow U_{f}$ is such that:
\begin{equation}
\underset{s\in \left[ 0,1\right] }{\sup }\int\nolimits_{0}^{T}\int%
\nolimits_{0}^{1}\left\| U_{f}\right\| ^{2}(r,s,t)r(1-r^{2})drdt\leq \dfrac{e%
}{4}\int\nolimits_{0}^{T}\int\nolimits_{0}^{1}\left\| \dfrac{\partial U_{s}%
}{\partial z}\right\| ^{2}dzdt.  \label{fluid1}
\end{equation}
\end{enumerate}
\end{Prop}

Let :
\begin{equation*}
\begin{array}{l}
W_{f}=L^{2}\left( 0,T;L^{2}\left( 0,1;W_{0}\right) \right) \cap L^{2}\left(
0,T;C\left( 0,1;\left( L^{2}\left( 0,1\right) \right) ^{N}\right) \right) ,
\\
W_{p}=L^{2}\left( 0,T;H^{1}\left( 0,1\right) \right) \cap C\left( 0,T;\left(
L^{2}\left( 0,1\right) \right) ^{N}\right) ,
\end{array}
\end{equation*}
and
\begin{equation*}
W_{g}=\left\{ \left( u,v\right) \in W_{f}\times W_{p}\mid u\left(
1,z,t\right) =v\left( z,t\right) \right\} .
\end{equation*}

\begin{Theo}
Under the hypothesis $\left( H3\right) $\ the problem (\ref{sys1}) admits a
solution if:
\begin{equation*}
\dfrac{\sqrt{e}}{2}\underset{i}{\sup }\left( \dfrac{\gamma _{is}}{\beta _{if}%
}\right) ^{1/2}<\text{ }\underset{i}{\inf }\theta _{is}.
\end{equation*}

This solution belongs to $W_{g}$.
\end{Theo}

\begin{proof}
Using (\ref{phi}) and (\ref{fluid1}), we obtain:

\begin{equation*}
\left\| U_{f}\right\| _{f}\leq \sqrt{\dfrac{e\mu }{\alpha ^{2}\left( 4-\mu
\alpha ^{2}\right) }}\left\| \tilde{U}_{f}\right\| _{f}\text{ ; }\left\|
U_{f}\right\| _{f}^{2}:=\underset{s\in \left[ 0,1\right] }{\sup }%
\int\nolimits_{0}^{T}\int\nolimits_{0}^{1}\left\| U_{f}\right\|
^{2}(r,s,t)r(1-r^{2})drdt.
\end{equation*}

The $\alpha $ which minimizes the Lipschitz constant is given by $\alpha
^{2}=2/\mu $, which leads to $\left\| U_{f}\right\| _{f}\leq \mu \sqrt{e}%
\left\| \tilde{U}_{f}\right\| _{f}/2$.\ This proves that $\Gamma :\tilde{U}%
_{f}\mapsto U_{f}$ is a contraction if and only if $\mu <2/\sqrt{e}.$
\end{proof}

\section{\textbf{Uniqueness of the solution}}

\begin{Theo}
Assuming that the solution of (\ref{sys1}) is smooth enough, the system (\ref
{sys1}) has at most one solution.
\end{Theo}

\begin{proof}
We suppose the existence of two couples of solutions $\left(
C_{if}^{1},C_{is}^{1}\right) _{i=1,...,N}$ and $\left(
C_{if}^{2},C_{is}^{2}\right) _{i=1,...,N}$ of (\ref{sys1}), and define: $%
W_{if}=C_{if}^{1}-C_{if}^{2}$ ; $W_{is}=C_{is}^{1}-C_{is}^{2}$. We multiply
the $i$-th equation of the system verified by $W_{if}$ by $rW_{if}$,
integrate on $\left[ 0,1\right] \times \left[ 0,1\right] \times \left[ 0,T%
\right] $ and use the equation verified by $W_{is}$ on the boundary. Thanks
to $\left( H3\right) $, one can deduce that $W_{if}$ is equal to zero in the
cylinder because $W_{if}$ is equal to zero at the inlet ($z=0$) and at the
outlet of the cylinder ($z=1$) and its partial derivative with respect to $r$
is equal to zero too, which implies, thanks to the equation verified by $%
W_{is}$ on the boundary, that the partial derivative of $W_{if}$ with
respect to $z$ is equal to zero in the cylinder. We also deduce that $\frac{%
\partial W_{is}}{\partial z}$ is equal to zero on the boundary. $W_{is}$\ is
equal to zero at time $T$ and at time $0$ because of the initial conditions.
Because of $T$ is arbitrarily chosen, $W_{is}$\ is equal to zero on the
boundary $\Sigma $ at any time.
\end{proof}

\section{Qualitative properties of the solution}

\subsection{\textbf{Nonnegativity of the solution}}

\begin{Prop}
For almost every $\left( r,z,t\right) $ in $\left] 0,1\right[ \times \left]
0,1\right[ \times \left] 0,T\right[ $, and for $i\in \left\{ 1,\ldots
,N\right\} $, we have : $0\leq C_{if}\left( r,z,t\right) $, $0\leq
C_{is}\left( z,t\right) $.
\end{Prop}

\begin{proof}
This is obtained multiplying the equations of (\ref{sys1}) by the
non-negative parts of $C_{if}$ or $C_{is}$, respectively.
\end{proof}

\subsection{\textbf{Upper and lower bounds of the concentrations}}

\begin{Prop}
\begin{enumerate}
\item  Let $\delta _{i}=-1$. For almost every $\left( r,z,t\right) $ in $%
\left] 0,1\right[ \times \left] 0,1\right[ \times \left] 0,T\right[ $, we
have : $0\leq C_{if}\left( r,z,t\right) \leq A_{i0}$ ; $0\leq C_{is}\left(
z,t\right) \leq A_{i0}$, with:
\begin{equation*}
A_{i0}=\max \left( \underset{r\in \left[ 0,1\right] }{\sup }C_{i0}\left(
r\right) ,\underset{z\in \left[ 0,1\right] }{\sup }C_{is0}\left( z\right)
\right) .
\end{equation*}

\item  Let $\delta _{i}=1$. For almost every $\left( r,z,t\right) $ in $%
\left] 0,1\right[ \times \left] 0,1\right[ \times \left] 0,T\right[ ,$ we
have: $a_{i0}\leq C_{if}\left( r,z,t\right) $ , $a_{i0}\leq C_{is}\left(
z,t\right) $, with:
\begin{equation*}
a_{i0}=\min \left( \underset{r\in \left[ 0,1\right] }{\inf }C_{i0}(r),%
\underset{z\in \left[ 0,1\right] }{\inf }C_{is0}(z)\right) .
\end{equation*}

\item  Let $\delta _{i}=1$. We have for almost every $\left( r,z,t\right) $
in $\left] 0,1\right] \times \left] 0,1\right[ \times \left] 0,T\right[ $ : $%
C_{if}\left( r,z,t\right) \leq a_{i0}e^{\lambda t}$ , $C_{is}\left(
z,t\right) \leq a_{i0}e^{\lambda t}$, with:
\begin{equation*}
a_{i0}=\max \left( \underset{r\in \left[ 0,1\right] }{\sup }C_{i0}(r),%
\underset{z\in \left[ 0,1\right] }{\sup }C_{is0}(z)\right) \text{ , }\lambda
=\sup_{i}k_{i},
\end{equation*}
$k_{i}$ being the Lipschitz constant of the function $\mathbf{r}_{i}$.

\item  Let $\delta _{i}=1$.\ There exist two positive constants $a$ and $b$
such that we have \ for every $l$ in $\left] 0,1\right[ $ and for every $T>0$%
:
\begin{equation*}
\begin{array}{rrl}
\int\nolimits_{0}^{T}\int\nolimits_{0}^{1}\left( C_{if}\right) ^{2}\left(
r,l,t\right) r\left( 1-r^{2}\right) drdt & \leq  & aT+b, \\
\int\nolimits_{0}^{1}\left( C_{is}\right) ^{2}\left( z,T\right) dz & \leq
& aT+b.
\end{array}
\end{equation*}
\end{enumerate}
\end{Prop}

\begin{proof}
The verification of these qualitative properties of the solution is
essentially obtained multiplying the equations of (\ref{sys1}) by the
appropriate non-negative or non-positive parts of the corresponding
test-functions.
\end{proof}

\section{\textbf{Numerical simulation for the reaction }$CO+O_{2}\rightarrow
CO_{2}$}

In the fluid, we use the following discretization method:

\begin{itemize}
\item  If $i$ is different of $N$ (chemical species), we directly solve (\ref
{sys1})$_{1}$ using the method of finite differences or that based on finite
elements.\ This requires that the equation with $i=N$ has already been
solved in order to put the appropriate values of $u_{N}$.

\item  If $i$ is equal to $N$ (temperature) we evaluate the coefficients at
the step before. This possibly requires the use of some fixed point argument.
\end{itemize}

On the boundary, we still use some finite differences method or some finite
element method (see \cite{crespim} for the details).

\medskip Let us take the following initial and boundary conditions:
\begin{equation*}
\begin{array}{cc}
\begin{array}{rrcl}
\text{in the cylinder:} & CO\left( r,0,t\right) & = & 0.02 \\
& O_{2}\left( r,0,t\right) & = & 0.05 \\
& CO_{2}\left( r,0,t\right) & = & 0 \\
& T\left( r,0,t\right) & = & 500,
\end{array}
&
\begin{array}{rrcl}
\text{on the boundary:} & CO\left( z,0\right) & = & 0.02 \\
& O_{2}\left( z,0\right) & = & 0.05 \\
& CO_{2}\left( z,0\right) & = & 0 \\
& T\left( z,0\right) & = & 490,
\end{array}
\end{array}
\end{equation*}

We have the following graphs at $t=0.3s$, $12s$, $24s$, $36s$ and $60s$.\ We
observe that the $CO$ and $O_{2}$ concentrations are decreasing and that the
temperature and the $CO_{2}$ concentration are increasing:

\bigskip

- Temperature 

\includegraphics[scale=.5]{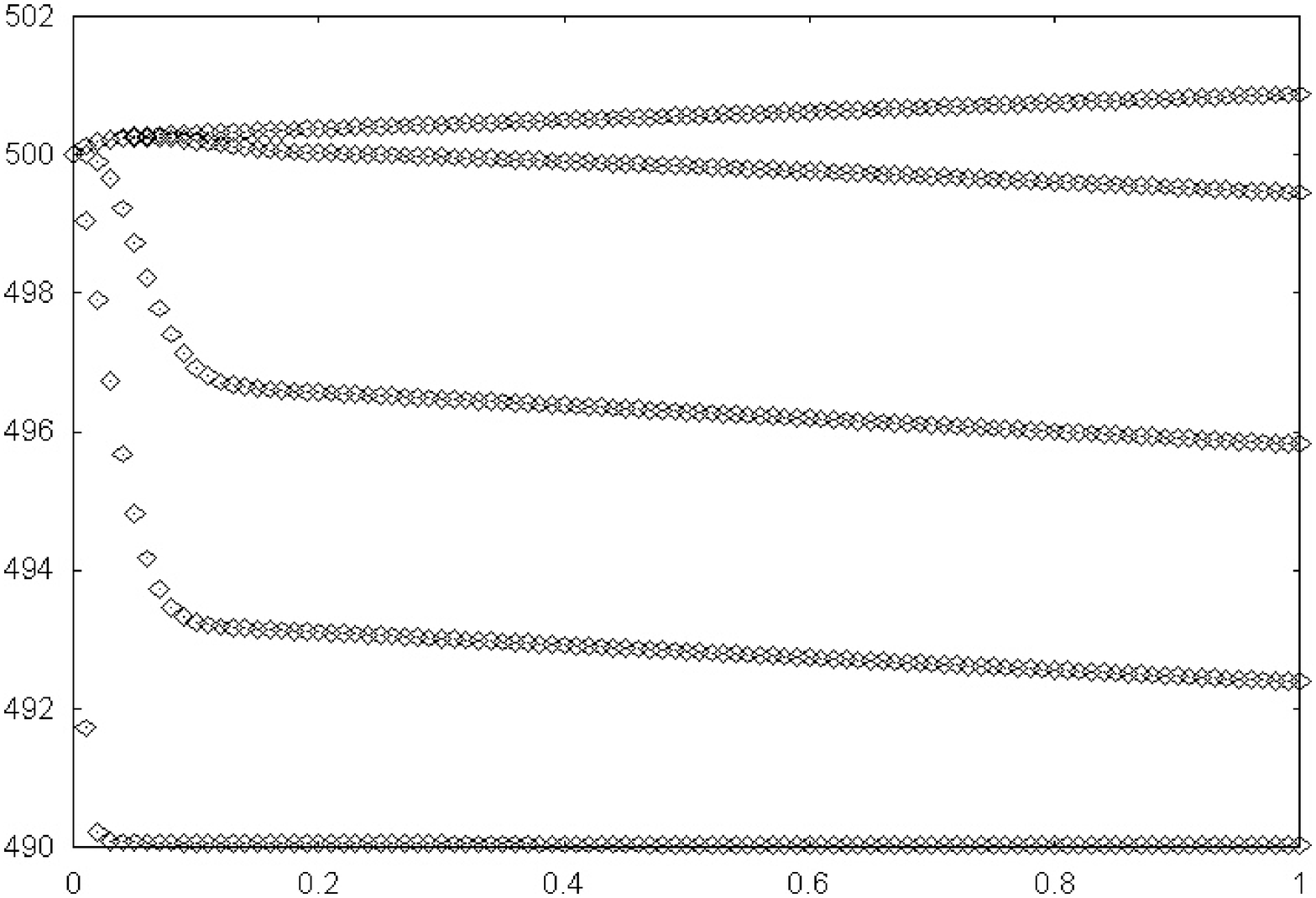}

\newpage

\textbf{- }$CO$

\includegraphics[scale=.5]{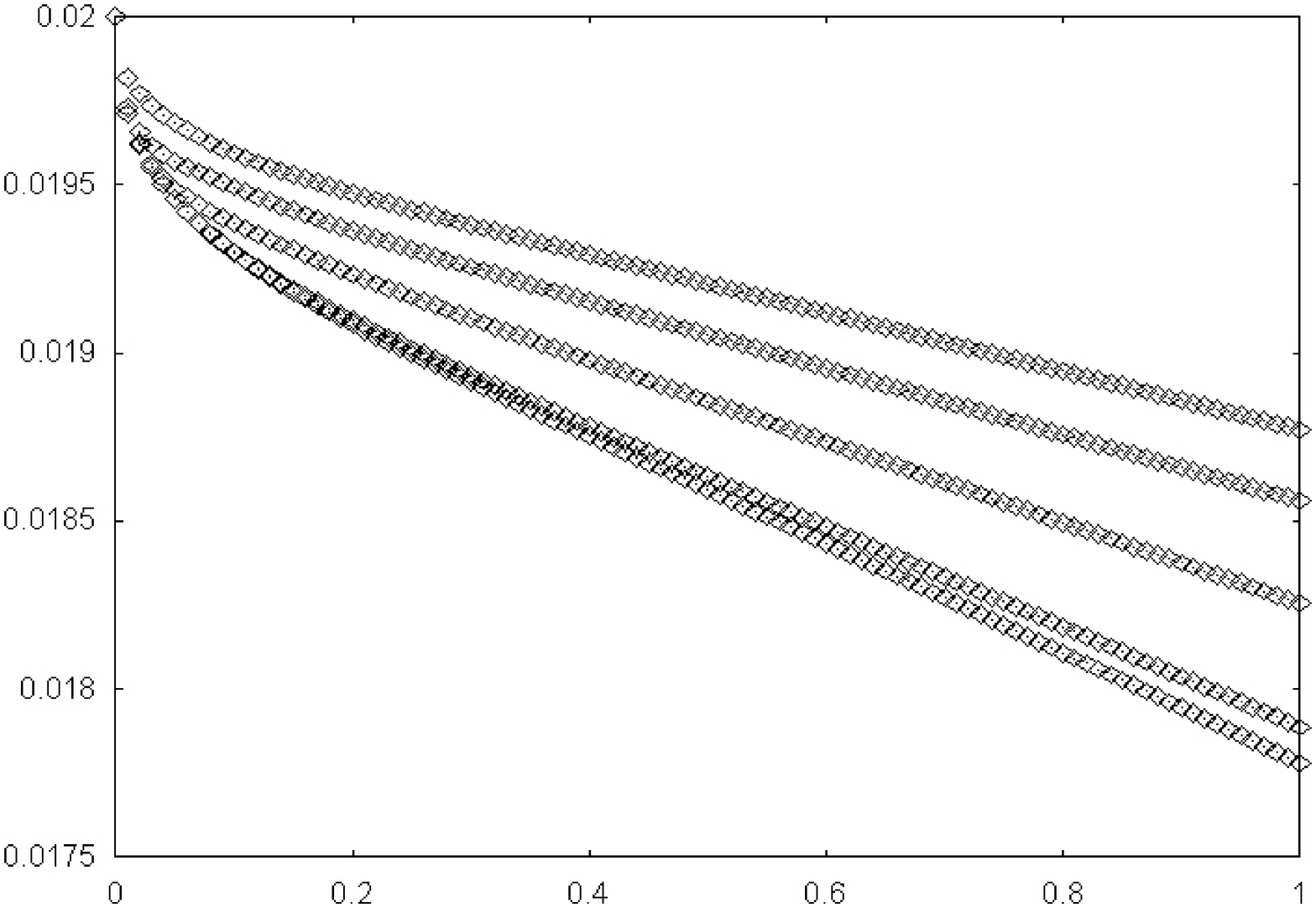}

- $O_{2}$

\includegraphics[scale=.5]{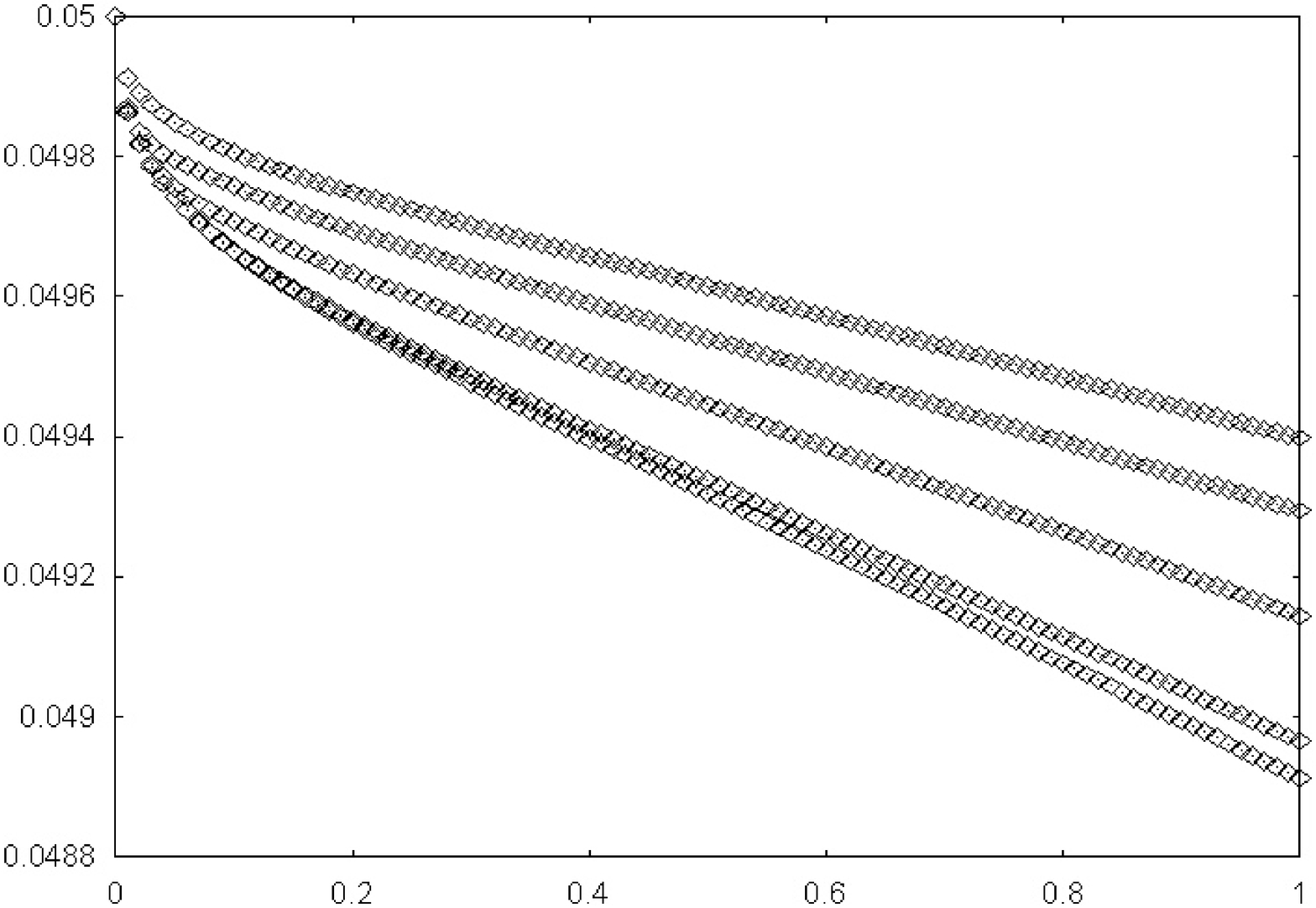}

\newpage
- $CO_{2}$

\includegraphics[scale=.5]{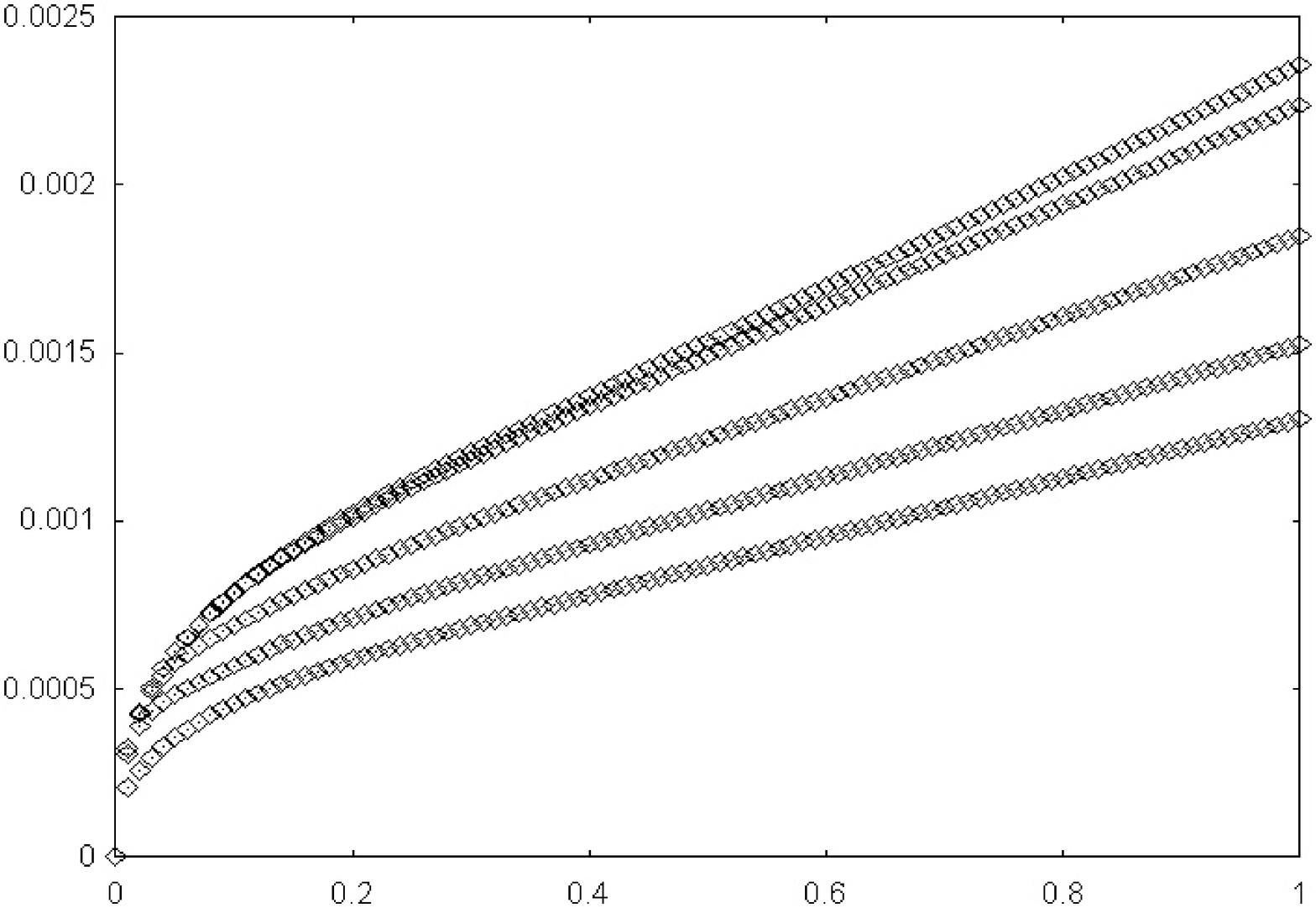}
\begin{Rem}
The reaction ends after 54s with the following values at the outlet of the
cylinder $(z=1)$:
\begin{equation*}
\begin{array}{rcl}
CO\left( 1,54\right) & = & 0.017777, \\
O_{2}\left( 1,54\right) & = & 0.048912, \\
CO_{2}\left( 1,54\right) & = & 0.002355, \\
T\left( 1,54\right) & = & 500.873497.
\end{array}
\end{equation*}
\end{Rem}

\end{document}